\documentclass[12pt]{article}
\usepackage{subfigure}
\usepackage{latexsym}
\usepackage{amsmath}
\usepackage{amssymb}
\usepackage{amsthm}
\usepackage{undertilde}
\usepackage{graphicx}
\usepackage{hyperref}
\usepackage[latin1]{inputenc}
\usepackage[T1]{fontenc}
\usepackage[all]{xy}
\newtheorem{lemma}{Lemma}
\newtheorem{example}{Example}
\newtheorem{theorem}{Theorem}
\newtheorem{corollary}{Corollary}
\newtheorem{definition}{Definition}
\newtheorem{prop}{Proposition}
\newtheorem{remark}{Remark}
\newcommand{\Res}{\mbox{Res}}
\makeatletter
\renewcommand{\@seccntformat}[1]
{\csname the#1\endcsname.\enspace} \makeatother
\def \XXint#1#2#3{{\setbox0=\hbox{$#1{#2#3}{\int}$}
     \vcenter{\hbox{$#2#3$}}\kern-.5\wd0}}

\begin{document}
\begin{center}
   {\bf On Approximating Ruin Probability of Double Stochastic Compound Poisson Processes\footnote{\today}}\\
{\sc Amir T. Payandeh Najafabadi$^{a,}$\footnote{Corresponding
author: amirtpayandeh@sbu.ac.ir; Phone No. +98-21-29903011; Fax No. +98-21-22431649}} \& Dan Kucerovsky$^{b}$\\
a Department of Mathematical Sciences, Shahid Beheshti University,
G.C. Evin, 1983963113, Tehran, Iran.\\
b Department of Mathematics and Statistics, University of New
Brunswick, Fredericton, N.B. CANADA E3B 5A3.
\end{center}
\begin{center}
{\sc Abstract}\end{center}\small Consider a surplus process which
both of collected premium and payed claim size are two independent
compound Poisson processes. This article derives two approximated
formulas for the ruin probability of such surplus process, say
double stochastic compound poisson process. More precisely, it
provides two mixture exponential approximations for ruin
probability of such double stochastic compound poisson process.
Applications to long--term Bonus--Malus systems and a heavy-tiled
claim size distribution have been given. Improvement of our
findings compare to the Cram\'er-Lundberg
upper bound has been given.\\
{\bf keywords:} Ruin probability; Double Stochastic Compound
Poisson Processes; Bonus--Malus system; Heavy-tailed
distributions; Laplace transforms; Cram\'er-Lundberg upper bound.
 \normalsize
\section{Introduction}
Consider double stochastic compound Poisson process
\begin{eqnarray} \label{surplus-process-BMS}
  U_t &=& u+\sum_{i=1}^{N_1(t)}C_i-\sum_{j=1}^{N_2(t)}X_j,
\end{eqnarray}
where  $C_1,C_2,\cdots$ and $X_1,X_2,\cdots,$ respectively, are
two independent i.i.d. random samples from independent random
premium $C$ and random claim size $X$, two independent Poisson
processes $N_1(t)$ and $N_2(t)$ (with intensity rates $\lambda_1$
and $\lambda_2$) are, respectively, stand for claims and purchase
request processes, and $u$ represents initial wealth/reserve $u$
of the process. Moreover, suppose that non-negative and continuous
random premium $C$ and claim $X,$ respectively, have density
functions $f_C$ and $f_X$ with some additional properties to be
discussed later (see Assumption 1 and Assumption 2).

The ruin probability for such process can be defined by
\begin{eqnarray} \label{Ruin-Definition}
  \psi(u) &=& p(T_{u}<\infty|U_{0}=u),
\end{eqnarray}
where $T_{u}$ is the hitting time, i.e., $T_{u}:=\inf\{t:~U_t\leq
0\}.$

Several authors studied the ruin probability under a stochastic
income assumption such as the surplus process
\ref{surplus-process-BMS}. For instance, Lappo (2004) provided an
equation for the ruin probability of a surplus process which
contains two compound Poisson processes. Rongming et al. (2007)
considered a surplus process with random premium and geometric
L\'evy investments return process. They obtained an
integro-differential equation for the ruin probability of such
process. Lim \& Qi (2009) considered a discrete-time surplus
process and established an equation for the ultimate ruin
probability. Moreover, they obtained an upper bound for the ruin
probability and studied its properties via a simulation study.
Under the compound binomial processes, Bao \& Wang (2013) provided
a difference equation and a defective renewal equation satisfied
respectively by the expected discounted penalty function. Another
version of the surplus compound binomial process has been studied
by Yu (2013a). Yu (2013b) considered surplus process that both
stochastic premium income and stochastic claims occurrence are
driven by the Markovian regime-switching process. Then, he derived
the Laplace transform of the time of ruin. Temnov (2014) provided
a recursive formula similar to the Beekman convolution formula to
evaluate the ruin probability for a stochastic premium surplus
process. Landriault \& Shi (2014) studied the finite-time ruin
probability for double stochastic compound poisson processes
whenever busy period follows a fluid flow model. Gatto \&
Baumgartner  (2014a) using the saddlepoint method to derive an
approximation for ruin probability of compound poisson risk
process perturbed by diffusion. Gatto \& Baumgartner (2014b)
extended Gatto \& Baumgartner (2014)'s findings to compound
poisson risk process perturbed by a Wiener process with infinite
time horizon. Cai \& Yang (2014) provided an integro-differential
equation for ruin probability of a compound Poisson surplus
process perturbed by diffusion with debit interest.

It is well-known that, in the situation that, random premium $C$
and random claim size $X$ are two independent exponential random
variables survival probability $\widetilde{\psi}(\cdot)$ can be
found explicitly as an exponential function (Melnikov, 2011,
Proposition 8.1). By an induction argument this fact can be
readily generalized to situation that $C$ and $X$ are two
independent mixture exponential random variables. In this
situation, survival probability $\widetilde{\psi}(\cdot)$ appears
as a mixture of some exponential functions. This article utilized
this fact and approximates survival probability
$\widetilde{\psi}(\cdot)$ appears as a mixture of some exponential
functions. The rest of this article is organized as follows. Some
mathematical background for the problem has been collected in
Section 2. Section 3 provides the main contribution of this
article. Applications of the results along with a and comparison
with the Cram\'er-Lundberg upper bound have been given in Section
4.
\section{Preliminaries}
The following recalls the exponential type $T$ functions which
plays a vital role in the rest of this article.
\begin{definition}
\label{exponential-type} An $L_1({\Bbb R})\cap L_2({\Bbb R})$
function $f$ is said to be of exponential type $T$ on ${\Bbb C}$
if there are positive constants $M$ and $T$ such that
$|f(\omega)|\leq M\exp\{T|\omega|\},$ for $\omega\in {\Bbb C}.$
\end{definition}
The Fourier transforms of exponential type functions are
continuous functions which are infinitely differentiable
everywhere and are given by a Taylor series expansion over every
compact interval, see Champeney (1987, page 77) and Walnut (2002,
page 81). The Fourier transform to exponential type functions are
also called band-limited functions, see Bracewell (2000, page 119)
for more details on band-limited functions. The well-known
Paley-Wiener theorem states that the Fourier transform of an
$L_2({\Bbb R})$ function vanishes outside of an interval $[-T,T],$
if and only if the function is an entire function of exponential
type $T$, see Dym \& McKean (1972, page 158) for more details. We
also need one form of the half-line version of the Paley-Wiener
theorem, namely that if a function in $f\in L_2({\Bbb R})$
vanishes on the left half-line, then its corresponding Fourier
transform $\hat{f}$ is holomorphic and uniformly bounded in the
upper half plane.  Ablowitz \& Fokas (1990 \S 4).

We also have need of the theory of residues for meromorphic
functions. Briefly, $\Res(f,z_i)$ is the coefficient of
$\frac{1}{z-z_i}$ in the Laurent series expansion of $f$ around
$z_i.$ It is possible to define Laurent series about the point
$z_i=\infty$, and hence it is possible to define a residue at
infinity, $\Res(f,\infty).$  The importance of the theory of
residues probably rests on two facts: many integrals can be
evaluated in terms of a sum of residues, and at poles of finite
order residues can be evaluated efficiently by formulas:
$$\Res(f,z_0)=\frac{1}{(m-1)!}\lim_{s\rightarrow z_0}\left(\frac{d}{ds}\right)^{m-1}\!\!\!\!\!\!\!\!\!(s-z_0)^m f(s),$$
where $m$ is the order of the pole, see for example Ablowitz \&
Fokas (1990 \S 4).

The following explores the Laplace transform of a function that
behaves like an exponential function about zero.
\begin{lemma}
\label{R_exponential} Suppose $f(\cdot)$ is an exponential type
function which behaves like an exponential function about zero,
i.e., $f(\epsilon)=\alpha_0e^{-\beta_0\epsilon},$ where $\alpha_0$
and $\beta_0$ are two given positive numbers. Then, the Laplace
transform $f(t+a)$ is
\begin{eqnarray}
\label{Laplace-f-t-s}
  \mathcal{L}(f(t+a);t;s) &=&
  e^{sa}\mathcal{L}(f(t);t;s)-\frac{\alpha_0}{s+\beta_0}\left(e^{as}-e^{-a\beta_0}\right).
\end{eqnarray}
\end{lemma}
\emph{Proof.} Suppose $f$ is of exponential type $T,$ the Taylor
series expansion for $f(t+a)$ is
$$f(t+a)=\sum_0^\infty  \frac{a^k}{k!} f^{(k)}(t),$$
where $f^{(k)}(\cdot)$ stands for the $k^{\hbox{th}}$ derivative
of $f(\cdot).$

Using the Weierstrass M-test along with properties of exponential
type function (i.e., $|f^{(k)}(t)|\leq |T|^k\|f\|_\infty.$). One
may conclude that the above Taylor series converges uniformly with
respect to both $a\in {\Bbb C}$ and $t\in{\Bbb R}.$ The uniform
convergence justifies the following taking term-by-term Laplace
transform.
\begin{eqnarray*}
  \mathcal{L} (f(t+a);t;s) &=& \mathcal{L}
  (f(t);t;s)+\sum_{k=1}^{\infty}\frac{a^k}{k!}\mathcal{L}
(f^{(k)}(t);t;s)\\
&=& \mathcal{L}
(f(t);t;s)\left[1+\sum_{k=1}^{\infty}\frac{(as)^k}{k!}\right]-\sum_{k=1}^{\infty}\frac{a^k}{k!}\sum_{j=0}^{k-1}s^jf^{(k-j-1)}(0)\\
&=& e^{as}\mathcal{L}
(f(t);t;s)-\sum_{k=1}^{\infty}\frac{a^k}{k!}\sum_{j=0}^{k-1}s^j\alpha_0(-\beta_0)^{k-j-1}\\
&=&e^{as}\mathcal{L}
(f(t);t;s)-\frac{\alpha_0}{s+\beta_0}\left(e^{as}-e^{-a\beta_0}\right),
\end{eqnarray*}
where the second equality arrived from an application of the
Laplace transform of a derivative, see Schiff (1999) for more
details. $\square$

We consider a class of random variables with their corresponding
density/probability functions and moment generating functions
satisfy the following two assumptions:
\begin{description}
    \item[Assumption 1.] Their density/probability and survival functions are exponential
        type functions, twice (continuously) differentiable, and bounded
           on ${\Bbb R}^+;$
    \item[Assumption 2.] Their moment generating functions are or can be extended to
    a meromorphic function on ${\Bbb C}.$
\end{description}
Exponential, mixture exponential, Erlang, Pareto (with finite
mean), normal, mixture normal, lognormal, etc are some
distributions that satisfy Assumption1 (Cai, 2004) while
Assumption 2 holds for most of the common distributions, such as
binomial, negative binomial, poisson, uniform, normal, chi-squared
(with an even degree of freedom), gamma (with integer shape
parameter), laplace, etc. On the other hand, several
nonmeromorphic moment generating functions can be analytically
continued to a meromorphic functions. Geometric and Multinomial
distributions are two examples which have nonmeromorphic moment
generating functions but their moment generating functions can be
analytically continued to a meromorphic function, see Kucerovsky
\& Payandeh (2014) and Sasv\a'ri (2013, \S 3.3) for more details.

By conditioning of the survival probability
$\widetilde{\psi}(\cdot)$ of surplus process
\ref{surplus-process-BMS} on the first arriving premium and claim
along with properties of compound Poisson process. One may derive
the following integro-differential equation for survival
probability $\widetilde{\psi}(\cdot)$ of surplus process
\ref{surplus-process-BMS}, proof may be found in Rongming et al.
(2007), Melnikov (2011, Theorem 8.1), among others.
\begin{equation}\label{integro-differential-equation}
-(\lambda_1+\lambda_2)\widetilde{\psi}(u)+\lambda_1
E(\widetilde{\psi}(u+C))+\lambda_2\int_0^{u}\widetilde{\psi}(u-x)f_X(x)dx=0,
\end{equation}
where $\lim_{u\rightarrow\infty}\widetilde{\psi}(u)=1$ and two
random premium $C$ and claim size $X$ satisfy assumptions
Assumption 1 and Assumption 2.

In the situation that, $C$ and $X$ are two independent exponential
random variables with rates $\mu_1$ and $\mu_2,$ respectively.
Then, survival probability can be found explicitly as
$\widetilde{\psi}(u)=Ae^{-Bu},$ where
$A:=(\mu_1+\mu_2)\lambda_2/(\lambda_1+\lambda_2)$ and
$B:=(\mu_1\lambda_1-\mu_2\lambda_2)/(\lambda_1+\lambda_2),$ see
Melnikov (2011, Proposition 8.1) for more details. This fact can
be generalized for situation that $C$ and $X$ are two independent
mixture exponential random variables.

Hereafter now, we assume for very small surplus/reserve, say
$u=\epsilon,$ survival probability behaves like an exponential
function, i.e.,
\begin{description}
    \item[Assumption 3.] For very small surplus/reserve $u=\epsilon$ survival probability  behaves like $\widetilde{\psi}(\epsilon)=\alpha_0e^{-\beta_0u},$
\end{description}
where two positive numbers $\alpha_0$ and $\beta_0$ have to be
determined. The above assumption can be explained by the above
observation along with the fact that a given density/probability
function can be approximated,  with some degree of accuracy, by a
mixture exponential distribution.

Coefficient $\beta_0$ in Assumption 3 can be found \emph{either}
by integro-differential Equation
\ref{integro-differential-equation} and letting
$u\rightarrow\infty$ \emph{or} the fact that
$\lim_{u\rightarrow\infty}\psi(u)e^{Ru}=c\in (0,1),$ where
adjustment coefficient $R$ is positive solution of
$\lambda_1M_C(R)+\lambda_2M_X(-R)=\lambda_1+\lambda_2$ (Kaas, et
al., Page 112), see Proposition \ref{prop:one_term_approx} for
more details. While coefficient $\alpha_0$ cannot be found using
integro-differential Equation \ref{integro-differential-equation}
or other theoretical assumptions on ruin probability. Therefore,
analogue to situation that $C$ and $X$ are two independent
exponential random variables, we set
$\alpha_0:=(1/E(C)+1/E(X))\lambda_2/(\lambda_1+\lambda_2).$
\section{Main results}
This section utilizes integro-differential Equation
\ref{integro-differential-equation} to derive two approximate
formulas for the ruin probability of a double stochastic compound
Poisson process. \ref{integro-differential-equation}. We seek an
analytical solution $\widetilde{\psi}(\cdot)$ which is an
exponential type function. In the other word, we assume:

\begin{description}
    \item[Assumption 4.] $|\widetilde{\psi}(\omega)|\leq Me^{T |\omega|},~\omega\in{\Bbb
    C},$ for some real numbers $M$ and $T$ in ${\Bbb R}.$
\end{description}
If this assumption is not met, as might be the case if, for
example, there are point masses in $\psi(\cdot),$ our methods  may
still result in a formal solution that can be verified by
substitution into integro-differential Equation
\ref{integro-differential-equation}.

The following from provides an explicit expression for the ruin
probability of a double stochastic compound Poisson process in an
analytical model.
\begin{theorem}
\label{New-Rongming-equation} Suppose $U_t$ represents double
stochastic compound Poisson process \ref{surplus-process-BMS}
which their random claim size $X$ and random premium $C$ satisfy
Assumption 1 and Assumption 2. Moreover, suppose that survival
probability $\widetilde{\psi}(\cdot)$ satisfies Assumption 3 and
Assumption 4. Then,
\begin{description}
    \item[i)] the Laplace transform of the ruin probability $\psi(\cdot)$
    satisfies \begin{eqnarray}\label{eqn:Rongming01}
  \mathcal{L}(\psi(u);u,s) &=&
  \frac{N(s)}{s(s+\beta_0)D(s)};
\end{eqnarray}
    \item[ii)] poles of $\frac{N(s)}{s(s+\beta_0)D(s)},$ say $z_j,$ can
be represented as $z_j:= -a_j \pm ib_j,$ where $a_j>0;$
\item[iii)] There is at most one simple pole on the negative real
axis,
\end{description}
where
$D(s):=-\lambda_1-\lambda_2+\lambda_1M_C(s)+\lambda_2M_X(-s),$
$N(s):=(s+\beta_0)D(s)-s\lambda_1\alpha_0(M_C(s)-M_C(-\beta_0)),$
and $M_C(\cdot)$ and $M_X(\cdot)$ represent moment generating
functions of $C$ and $X,$ respectively.
\end{theorem}
\emph{Proof.} Taking a Laplace transform from integro-differential
Equation \ref{integro-differential-equation} with an application
of Lemma \ref{R_exponential} complete part (i). For part (ii)
observe that: by Assumption 3 and the Paley-Wiener theorem, we may
conclude that an expression $g(s):=\frac{N(s)}{s(s+\beta_0)D(s)}$
is holomorphic and uniformly bounded in the right-half plane.
Therefore, $g(s)$ has no poles at any $z\in{\Bbb C}$ with
non-negative real part. On the other hand, since $\psi(\cdot)$ is
real-valued function, its corresponding Laplace transform is real
everywhere on the real line (except possibly at poles) and hence
by the Schwartz reflection principle,
 $$\overline{\left(\frac{N(s)}{s(s+\beta_0)D(s)}\right)}=\frac{N(\overline{s})}{\overline{s}(\overline{s}+\beta_0)D(\overline{s})}.$$
Therefore, poles of $g(s)$ are located at  $z_j:= -a_j \pm ib_j.$
To establish $a_j>0,$ observe that if there is a pole at some
$z_i\in{\Bbb C},$ then there is also a pole at the complex
conjugate, $\overline{z_i}.$ On the other hand, since $D(s)$ is a
convex (concave up) function on the real line, it can have at most
two zeros on the real line. Moreover, $D(s)$ has a zero at the
origin, so there can be at most one other zero of $D(s)$ on the
real line, necessarily located on the negative real line. For part
(iii) observe that a simple pole for $g(s)$ on the negative real
axis can be $-\beta_0.$ $\square$

Using the above theorem, one need only consider the roots of
$D(s)=0,$ which are completely determined by the moment generating
functions, random claim size $X,$ and random random premium $C,$
say respectively $M_X(\cdot)$ and $M_C(\cdot).$ Residues at simple
poles are particularly easy to evaluate, and the case of simple
poles is the generic case. If there are finitely many poles and
all are of finite order, they may be made to be simple poles by an
infinitesimal perturbation of the problem. Thus the following
result is of practical importance. The condition that appears
below, of the derivatives being non-zero, is more or less
equivalent to all the zeros of $D(s)$ being simple zeros:
 \begin{corollary}
 \label{psi_u}
Under the same conditions as in Theorem
\ref{New-Rongming-equation}, we then have that the ruin
probability $\psi(\cdot)$ is
\begin{eqnarray*}
\psi(u) &\approx& \sum_{D(z_i)=0,\, \Re(z_i)<0}
\frac{N(z_i)}{z_i(z_i+\beta_0)D^\prime(z_i)} e^{z_iu}.
\end{eqnarray*} whenever the derivatives that appear are all non-zero.
\end{corollary}
\emph{Proof.} The residue of $N(s)/(s(s+\beta_0)D(s))$ at a simple
zero $z_i$ of $D(s)$ is $$\lim_{s\rightarrow
z_i}\frac{(s-z_i)N(s)}{s(s+\beta_0)D(s)}=\frac{N(z_i)}{z_i(z_i+\beta_0)D^\prime(z_i)}.$$
Conversely, if $D'(s)$ is nonzero at a zero of $D(s)$ then this
zero is a simple zero, and $N(s)/(s(s+\beta_0)D(s))$  has at most
a simple pole (as the zero is assumed to have negative real part).
$\square$

Since the ruin probability is real and non-negative, and since the
poles come in pairs, we can draw some further conclusions about
the structure of the expression for $\psi(\cdot).$
\begin{prop}
\label{prop:approx.by.exp} Under the same conditions as in Theorem
\ref{New-Rongming-equation}. The ruin probability $\psi(u)$ can be
simplified by a finite sum of the form
\begin{eqnarray}
\label{Sin-cos-ruin} \psi(u) &\approx& \alpha_0 e^{-c_0 u}
+\sum_{i} A_i e^{-c_i u} \cos(d_i u) + \sum_{i}B_i e^{c_i u}
\sin(d_i u),
\end{eqnarray}
where $c_i$ and $d_i$ are, respectively, real and imaginary parts
of the roots of $D(s)=0$ and $A_i,$ $B_i$ are real numbers
evaluated by substituting \ref{Sin-cos-ruin} in Equation
\ref{integro-differential-equation}.
\end{prop}
\emph{Proof.} We have the above general form of a sum involving
trigonometric functions and exponentials comes directly from
Theorem \ref{New-Rongming-equation}, plus the fact that the poles
come in pairs (Theorem \ref{New-Rongming-equation}.) The $b_i$ are
the imaginary parts of the roots, and the $a_i$ are the real
parts. Theorem \ref{New-Rongming-equation} moreover tells us that
there is at most one root on the negative real axis. Since the
ruin probability must remain non-negative for large $u,$ but must
vanish at infinity, it follows that the exponential corresponding
to this purely real root must be the dominant term for large $u,$
and thus the purely real root is negative, nonzero, and is to the
right of all the complex roots.  $\square$
\begin{remark}
It is a corollary of the proof above  that under our conditions
there does always exist a unique negative and nonzero real root of
$D(s).$
\end{remark}
\begin{corollary}
\label{Constant-2-zero} Under the same conditions as in Theorem
\ref{New-Rongming-equation}, and also assuming $D(s)/s=0$ has only
a simple non-zero root at $z_1,$ then, the ruin probability is
given by
\begin{eqnarray*}
  \psi(u) &\approx& \frac{\lambda_2}{\lambda_1+\lambda_2-\lambda_1 M_C (z_1)}e^{z_1 u},
\end{eqnarray*}
where $M_C(z_1)$ stands for the moment generating function of the
random variable $C$ at $z_1.$
\end{corollary}
\emph{Proof.} Substituting $u=0$ in Equation
\ref{integro-differential-equation}, one may conclude that
\begin{eqnarray*}
  \alpha_0 &:=& \psi(0)\\
  &=&\frac{\lambda_1}{\lambda_1+\lambda_2}E(\psi(C))+\frac{\lambda_2}{\lambda_1+\lambda_2}\\
  &=&  \frac{\lambda_1}{\lambda_1+\lambda_2}\frac{N(z_1)}{-z_1D^\prime
  (z_1)}M_C(z_1)+\frac{\lambda_2}{\lambda_1+\lambda_2}\\
  &=& \frac{\lambda_1}{\lambda_1+\lambda_2}\alpha_0M_C(z_1)+\frac{\lambda_2}{\lambda_1+\lambda_2},
\end{eqnarray*}
where two last equalities arrive from double applications of
Corollary \ref{psi_u}, with simple non-zero root $z_1$ at $u=0.$
$\square$

In general, the moment generating function of a distribution, when
it exists, will grow rapidly as we move to the right on the real
axis, but decreases along the negative real axis, and decreases
along the imaginary axis (because there it is the characteristic
function). Thus, it is reasonable to expect that $M_C (z_k)$ will
be small for the off-axis roots of $D(s)$, since these roots are
located in the second and third quadrants of the complex plane.
Furthermore, since the moment generating function oscillates in
the off-axis direction, sums involving the off-axis roots are
likely to be quite small due to cancellations. If this is the
case, then we can hope to approximate by neglecting the sum
involving off-axis roots roots in the above calculation. We thus
obtain the following approximation result.
\begin{prop} \label{prop:one_term_approx}.
Under the same conditions as in Theorem
\ref{New-Rongming-equation}, then the ruin probability
$\psi(\cdot)$ about the origin can be approximated by one of the
following exponential functions.
\begin{description}
    \item[i)] If equation $D(s)=0$ has unique strictly
negative real root, say $z_1,$
\begin{eqnarray*}
  \psi(\epsilon) &\approx& \frac{\lambda_2}{\lambda_1+\lambda_2-\lambda_1 M_C (z_1)}e^{z_1
  \epsilon};
\end{eqnarray*}
    \item[ii)] If $\lambda_1E(C)\geq \lambda_2E(X)$ \begin{eqnarray*}
  \psi(\epsilon) &\approx&
  \frac{E(C)+E(X)}{E(C)E(X)(\lambda_1+\lambda_2)}e^{-\beta_0
  \epsilon};
\end{eqnarray*}
    \item[iii)]  \begin{eqnarray*}
  \psi(\epsilon) &\approx&
  \frac{E(C)+E(X)}{E(C)E(X)(\lambda_1+\lambda_2)}e^{-R
  \epsilon},
\end{eqnarray*}
\end{description}
where $\beta_0$ is positive solution of
$\lambda_1M_C(-\beta_0)+\lambda_2M_X(\beta_0)=\lambda_1+\lambda_2$
and adjustment coefficient $R$ is positive solution of
$\lambda_1M_C(R)+\lambda_2M_X(-R)=\lambda_1+\lambda_2.$
\end{prop}
\emph{Proof.} Part (i) arrives by an application of Corollary
\ref{Constant-2-zero}. For part (ii) observe that substituting an
exponential function $\alpha_0e^{-\beta_0u}$ into
integro-differential Equation \ref{integro-differential-equation}
and letting $u\rightarrow\infty.$ One may obtain equation
$\lambda_1M_C(-\beta_0)+\lambda_2M_X(\beta_0)=\lambda_1+\lambda_2.$
Now, set
$h(\beta_0):=\lambda_1M_C(-\beta_0)+\lambda_2M_X(\beta_0).$ Using
assumption on random premium $C$ and random claim size $X,$ one
may conclude that $h(0)=\lambda_1+\lambda_2,$
$\lim_{\beta_0\rightarrow\infty}h(\beta_0)=\infty,$
$\lim_{\beta_0\rightarrow\infty}h^\prime(\beta_0)=-\lambda_1E(C)+\lambda_2E(X)\leq0,$
and $h^{\prime\prime}(\beta_0)\geq0.$ Therefore, the above
equation has a unique positive solution in $\beta_0.$ Part (iii)
arrives by an application of the fact that
$\lim_{u\rightarrow\infty}\psi(u)e^{Ru}=c\in (0,1),$ (Kaas, et
al., Page 112). Coefficient $\alpha_0$ for both parts (ii \& iii)
cannot be found using integro-differential Equation
\ref{integro-differential-equation} or other theoretical
assumptions on ruin probability. Therefore, analogue to situation
that $C$ and $X$ are two independent exponential random variables,
we set $\alpha_0:=(1/E(C)+1/E(X))\lambda_2/(\lambda_1+\lambda_2).$
$\square$

Residues at simple poles are particularly easy to evaluate, and
the case of simple poles is the generic case. If there are
finitely many (finite order) poles. They may be made to be simple
poles by an infinitesimal perturbation of the problem. Therefore,
an approximation for the ruin probability can be obtained by
assuming all roots of $D(s)=0$ are simple. The following theorem
improves accuracy of approximation result given by Theorem
\ref{New-Rongming-equation} whenever $D(s)=0$ has finitely many
simple zeros.
\begin{theorem}
\label{linear_systems} Supposing that under the conditions of the
above result, $D(s)$ has finitely many roots in the left
half-plane.  Then, $\psi(u)\approx\sum A_i e^{z_i u},$ where $A_i$
can be found out by the following linear system of equations
\begin{equation}
\label{eq:linearsystem} A_j
=\frac{1}{\lambda_2z_jM^\prime_X(-z_j)}\left(\lambda_2M_X(-z_j)-\lambda_2+z_j\sum_{i\neq
j}\frac{A_i}{z_j-z_i}\left(\lambda_1+\lambda_2-\lambda_1M_C(z_i)-\lambda_2M_X(-z_j)\right)
\right).
\end{equation}
\end{theorem}
\emph{Proof.} Roots of $D(s)=0$ may be made to be simple poles by
an infinitesimal perturbation of the problem. Therefore,
$D(\cdot)$ has finitely many simple roots, $z_i$, we have by Lemma
\ref{psi_u} and Proposition \ref{prop:approx.by.exp} that
$\psi(u)\approx\sum A_i e^{z_i u}.$ Substituting
$\psi(u)\approx\sum A_i e^{z_i u}$ into Equation
\ref{integro-differential-equation} and taking the Laplace
transform, we obtain
\begin{eqnarray*}
-(\lambda_1+\lambda_2)\left(\frac{1}{s}-\sum_{i}\frac{A_i}{z_j-z_i}
\right)+\lambda_1\left(\frac{1}{s}-\sum_{i}\frac{A_i}{z_j-z_i}M_C(z_i)
\right)+\lambda_2\left(\frac{1}{s}-\sum_{i}\frac{A_i}{z_j-z_i}
\right)M_X(-s) &=&  0.
\end{eqnarray*}
The desired proof arrives by setting $s=z_j$ along with an
application of the $l'H\hat{o}pital'$s rule. $\square$

It is clear from the above that the solution obtained is in fact
unique. The only further point to mention with respect to the
above theorem that since we now are dealing with approximations,
it could happen that the sum $\sum A_i e^{z_i u}$ obtained as in
Theorem \ref{linear_systems} is not quite real or not quite
non-negative, so we adjust it in the above by taking the real part
and then the floor with zero. The amount by which the sum $\sum
A_i e^{z_i u}$ fails to be real and non-negative can be taken as a
practical indication of the amount of error in the approximation.
Normally, one would pick the complex zeros that are closest to the
real root of $D(s),$ aiming to ensure that the complex moment
function $M_C(z)$ has quite small absolute value at the zeros that
have been neglected. Moreover, the set of roots chosen should be
closed under complex conjugation, meaning that if $a+ib$ is in the
set of roots, then $a-ib$ should also be in the set of roots.

It would be worthwhile mentioning that, the ruin probability may
be impacted by different ways of adjustment (i.e., adjusting {\it
either} $\psi(u)$ {\it or} complex roots). Sensitivity of such
adjustment should be investigated and managed by the amount of
error which obtained by substituting adjusted $\psi(\cdot)$ into
Equation \ref{integro-differential-equation}.

The following section is devoted to applications of Corollary
\ref{Constant-2-zero} and Theorem \ref{linear_systems}.
\section{Applications}
In the most of practical application the ruin probability of
surplus process \ref{surplus-process-BMS} has been approximated by
the Cram\'er-Lundberg upper bound $e^{-Ru},$ where adjustment
coefficient $R$ is positive solution of
$\lambda_1M_C(R)+\lambda_2M_X(-R)=\lambda_1+\lambda_2.$ Certainly,
$-R$ is one of solution of our equation $D(s)=0.$ Practical
implementation of the above findings along with a comparison with
the well-known Cram\'er-Lundberg upper bound have been given in
the following examples. Example 1 considers a situation that
equation $D(s)=0$ has just one simple non-positive solution. In
this situation our method improves the Cram\'er-Lundberg upper
bound by a constant coefficient. While in the Example 2 that
equation $D(s)=0$ has more than one simple non-positive root, our
method provides a significant improvement on the Cram\'er-Lundberg
upper bound.
\begin{example}
Suppose an insurance company plans to charge its policyholders
based upon their risks (such insurance systems well known as
Bonus--Malus systems). Moreover suppose that the insurance company
considers three different scenarios with 5, 10, and 20 premiums'
values such that under all of these scenarios expected revenue of
the company stay the same by selling an insurance contract. Based
upon a statistical investigation an actuary suggested the
following premium values and probability that a given policyholder
falls in a given level, say ${\bf \pi},$ of each
scenario\footnote{These three different scenarios also can be
viewed as three different Bonus--Malus systems which stabilized,
in the long run, around their corresponding equilibrium
distributions ${\bf \pi},$ well known as long--term Bonus--Malus
systems.}.
\begin{center}
\tiny Table 1: Premiums values ${\bf C}$ and probability that
a given policyholder falls in a given level for such three scenarios.\normalsize\\
\begin{tabular}{ c c c}
  \hline
Scenario (number of its levels) &  ${\bf \pi}$ & Premium value ${\bf C}$ \\
  \hline
S1  (with 5 levels) &  $\pi_i\equiv0.2$ &(0.6, 1, 1.4, 1.8, 2.2)\\
S2  (with 10 levels) & $\pi_i\equiv0.1$ & (0.4, 0.8, 1, 1.2, 1.4, 1.5, 1.7, 1.8, 2, 2.2) \\
S3  (with 20 levels) & $\pi_i\equiv0.05$ & (0.4, 0.5, 0.6, 0.7, 0.8, 0.9, 1, 1.1, 1.2, 1.3,\\
  --&  --           & 1.4, 1.5, 1.6, 1.7, 1.8, 1.9, 2.1, 2.2, 2.6, 2.7) \\
  \hline
\end{tabular}
\end{center}
Moreover, suppose that the number of sold contracts, $N_1(t),$ and
number of arrived claims, $N_2(t),$ are two independent Poisson
processes with intensity $\lambda_1=18$ and $\lambda_2=11$ and
claim size distributions given by Table 2.

It is easy to show that, for all cases, equation $D(s)=0$ has just
one simple non-positive root. Using results of Corollary
\ref{Constant-2-zero}, one may evaluate the ruin probability in an
exact mode. The ruin probability (as well as the Cram\'er-Lundberg
upper bound) for such three different scenarios have been given by
Table 2.
\begin{center}
\tiny Table 2: The ruin probability and Cram\'er-Lundberg upper
bound for the above three scenarios
under different claim size distributions.\normalsize\\
\begin{tabular}{c c c c}
  \hline
  & & Scenario  (number of its levels)\\
\cline{2-4}
Claim size distribution & Scenario 1 (5) &Scenario 2 (10) & Scenario 3 (20)\\
  \hline
  Gamma(1,3)  & $\psi(u)=0.313e^{-1.760u}$ & $\psi(u)=0.312e^{-1.759u}$ & $\psi(u)=0.280e^{-1.746u}$ \\
              & $(\psi(u)\leq e^{-1.760u})$ & $(\psi(u)\leq e^{-1.759u})$ & $(\psi(u)\leq e^{-1.746u})$\\
  Gamma(1,1)  & $\psi(u)=0.072e^{-0.408u}$ & $\psi(u)=0.072e^{-0.407u}$ & $\psi(u)=0.071e^{-0.402u}$ \\
              & $(\psi(u)\leq e^{-0.408u})$ & $(\psi(u)\leq e^{-0.407u})$ & $(\psi(u)\leq e^{-0.402u})$ \\
  Gamma(3,2)  & $\psi(u)=0.057e^{-0.229u}$ & $\psi(u)=0.057e^{-0.227u}$ & $\psi(u)=0.056e^{-0.223u}$ \\
              & $(\psi(u)\leq e^{-0.229u})$ & $(\psi(u)\leq e^{-0.227u})$ & $(\psi(u)\leq e^{-0.223u})$ \\
  Gamma(5,3)  & $\psi(u)=0.053e^{-0.174u}$ & $\psi(u)=0.053e^{-0.173u}$ & $\psi(u)=0.053e^{-0.170u}$ \\
              & $(\psi(u)\leq e^{-0.174u})$ & $(\psi(u)\leq e^{-0.173u})$ & $(\psi(u)\leq e^{-0.170u})$ \\
  \hline
\end{tabular}
\end{center}
As one may observe, our method just improve the Cram\'er-Lundberg
upper bound by a constant coefficient.

Figure 1 illustrates several comparison regarding to the ruin
probability of such three scenarios under different Claim size
distributions.
\begin{center}
\begin{figure}[h!]
\centering\subfigure[]{
\includegraphics[width=4cm,height=5cm]{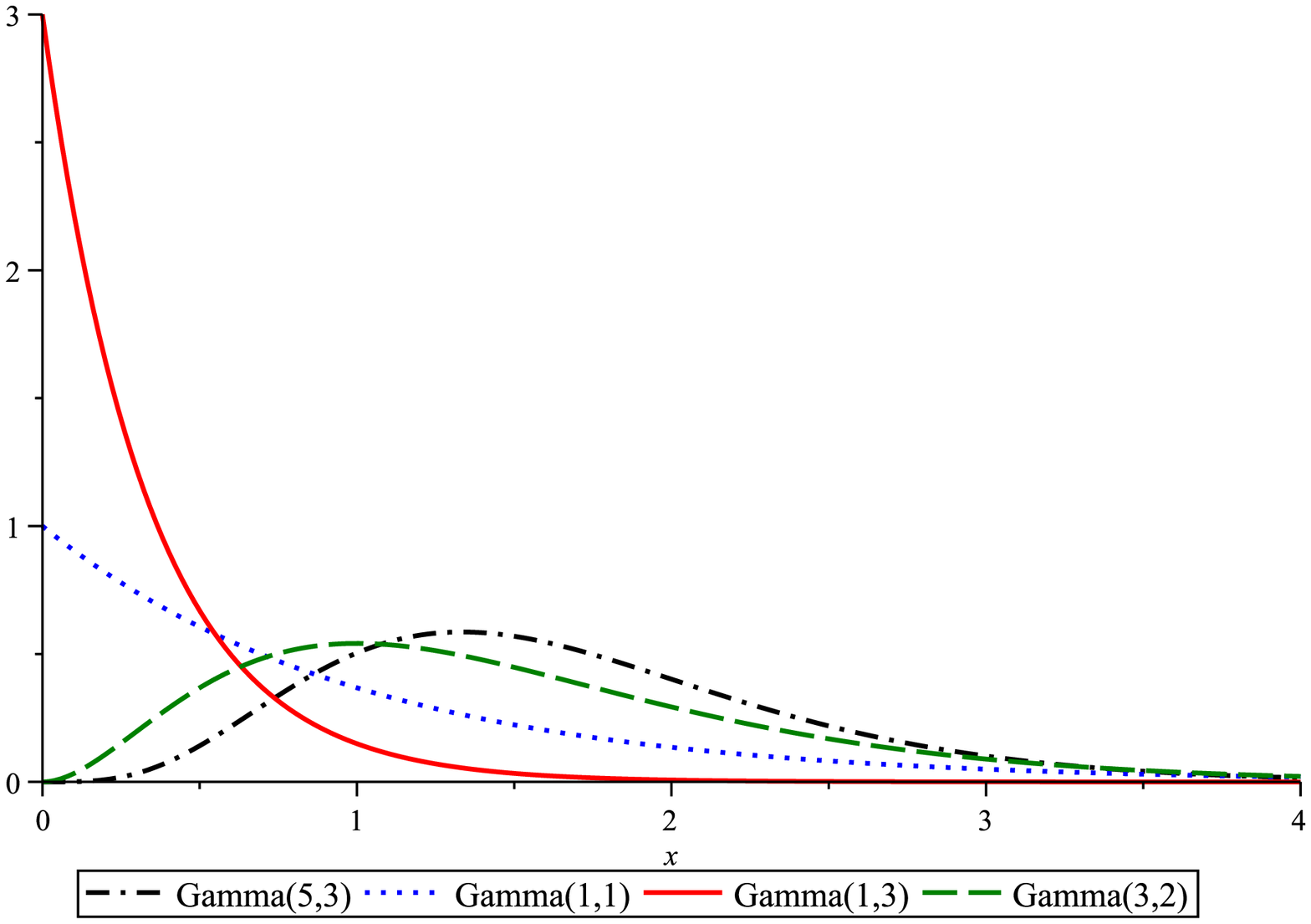}}\subfigure[]{
\includegraphics[width=4cm,height=5cm]{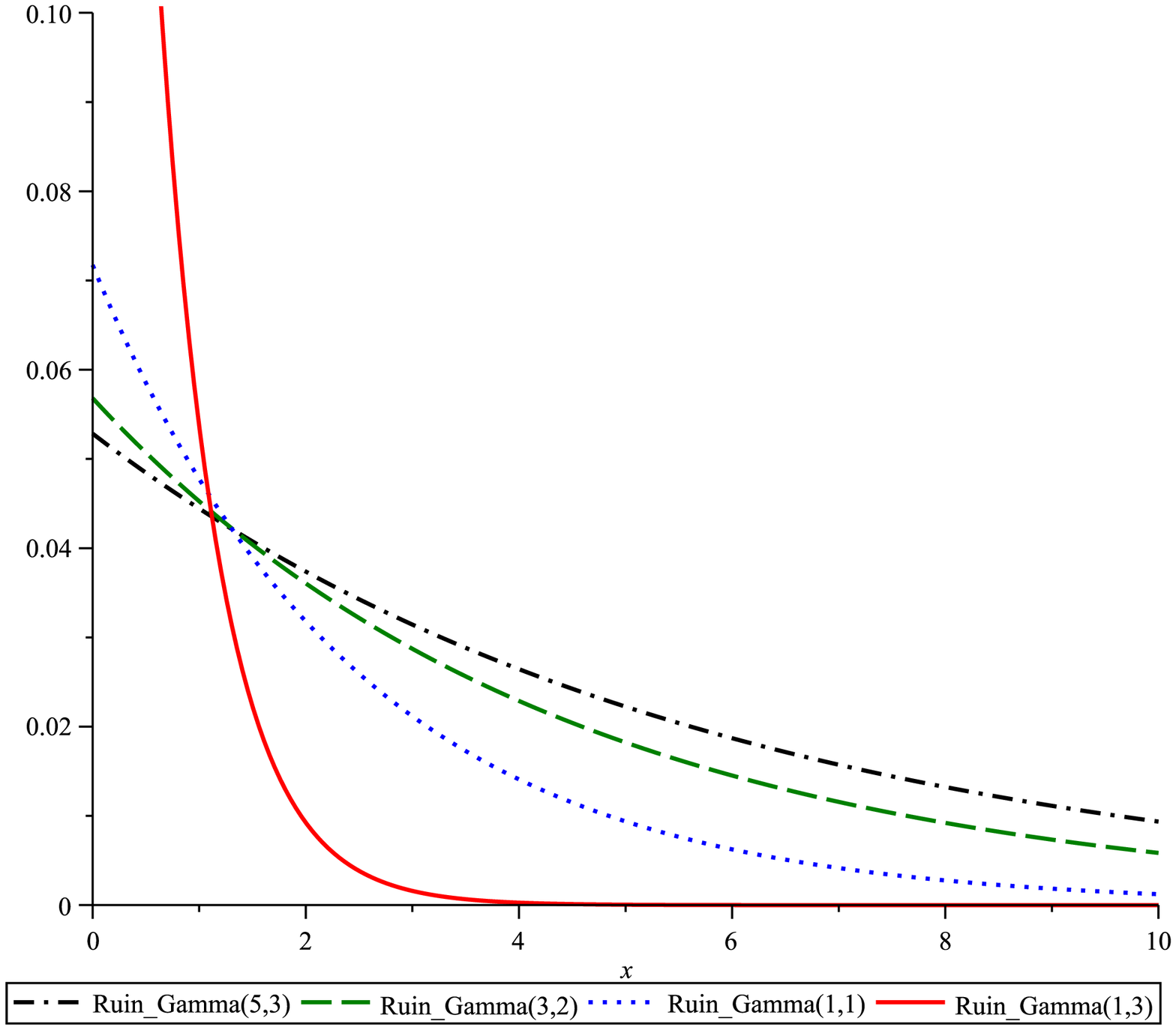}}\subfigure[]{
\includegraphics[width=4cm,height=5cm]{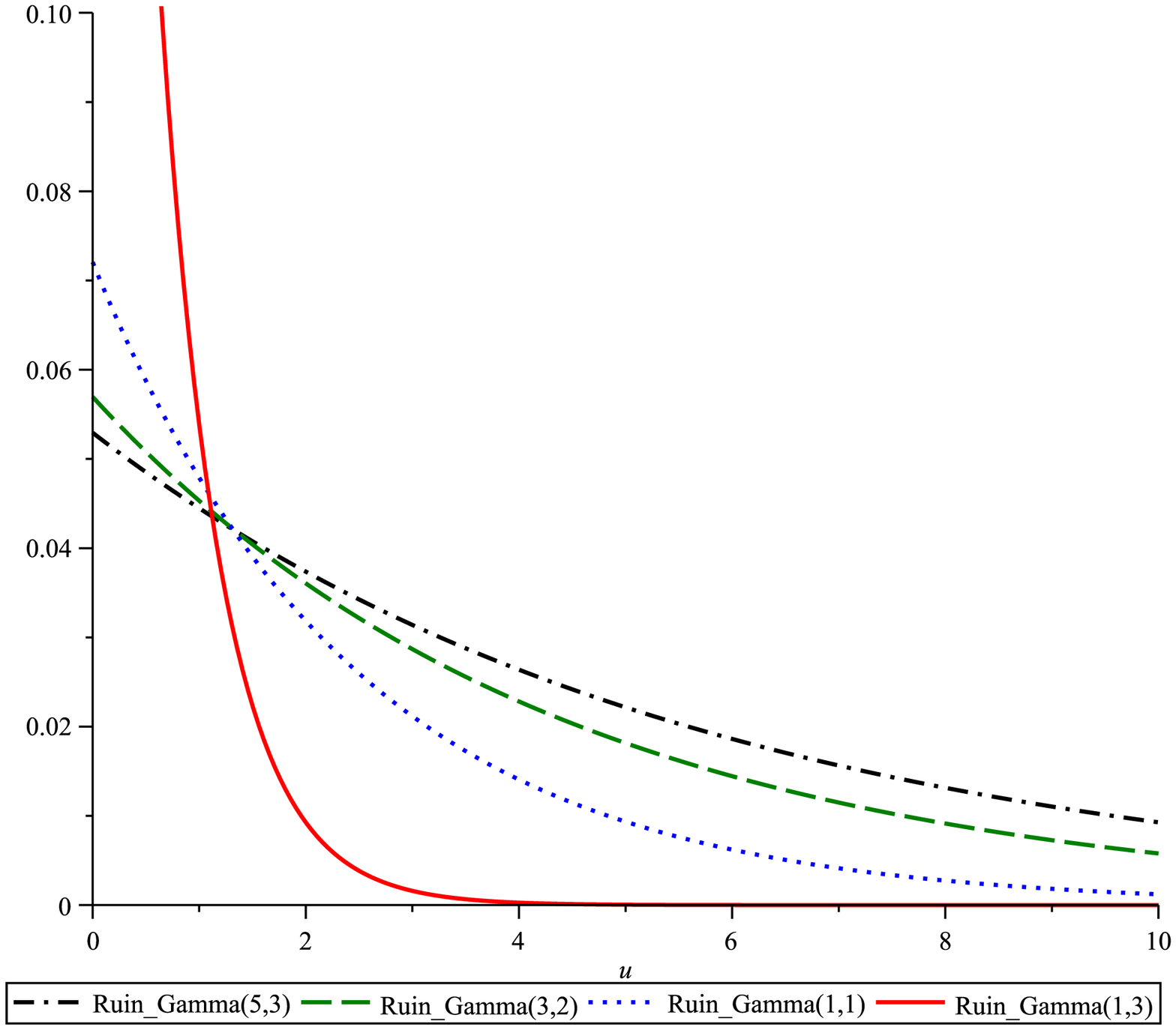}}\subfigure[]{
\includegraphics[width=4cm,height=5cm]{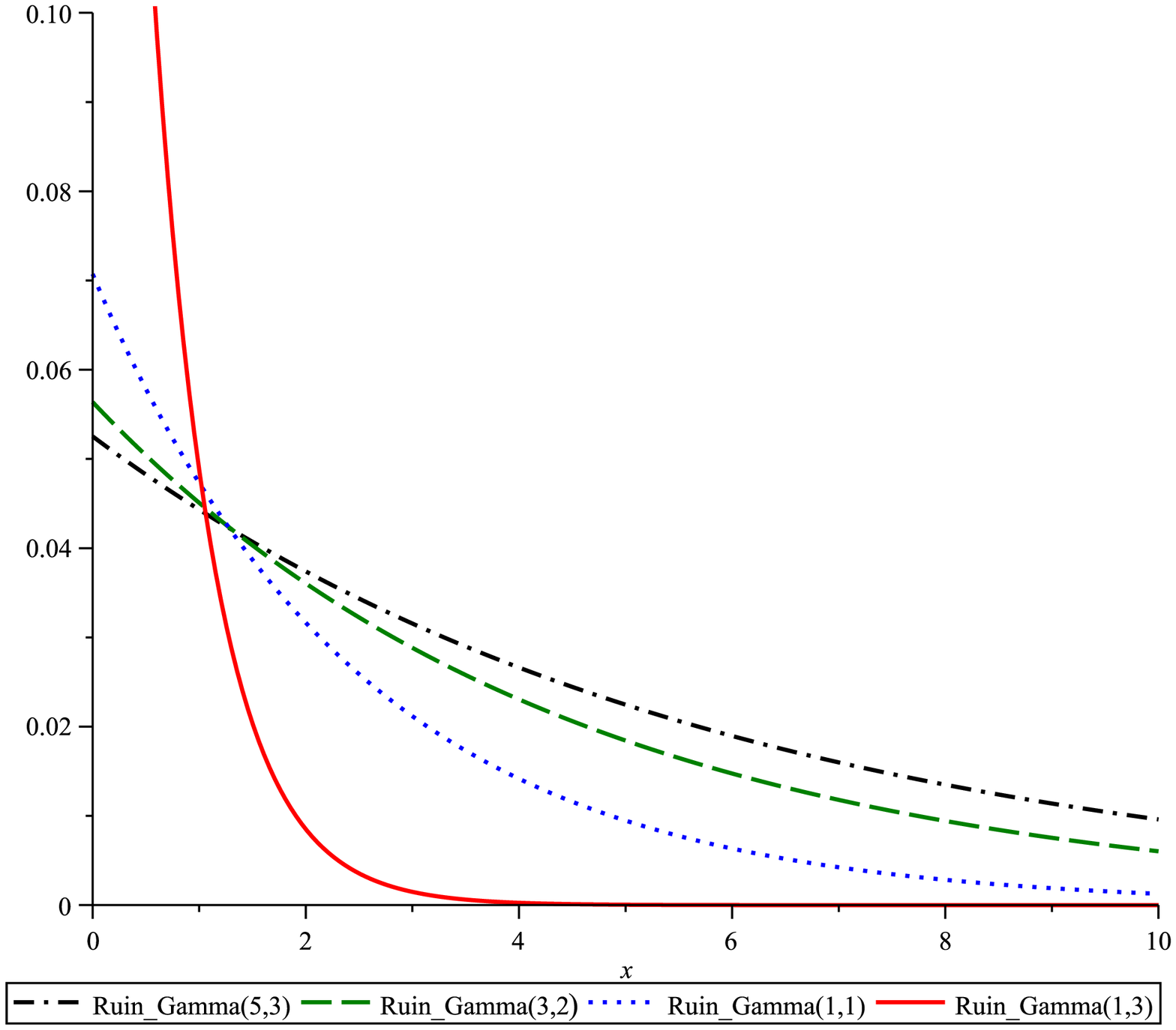}}
\centering\subfigure[]{
\includegraphics[width=4cm,height=5cm]{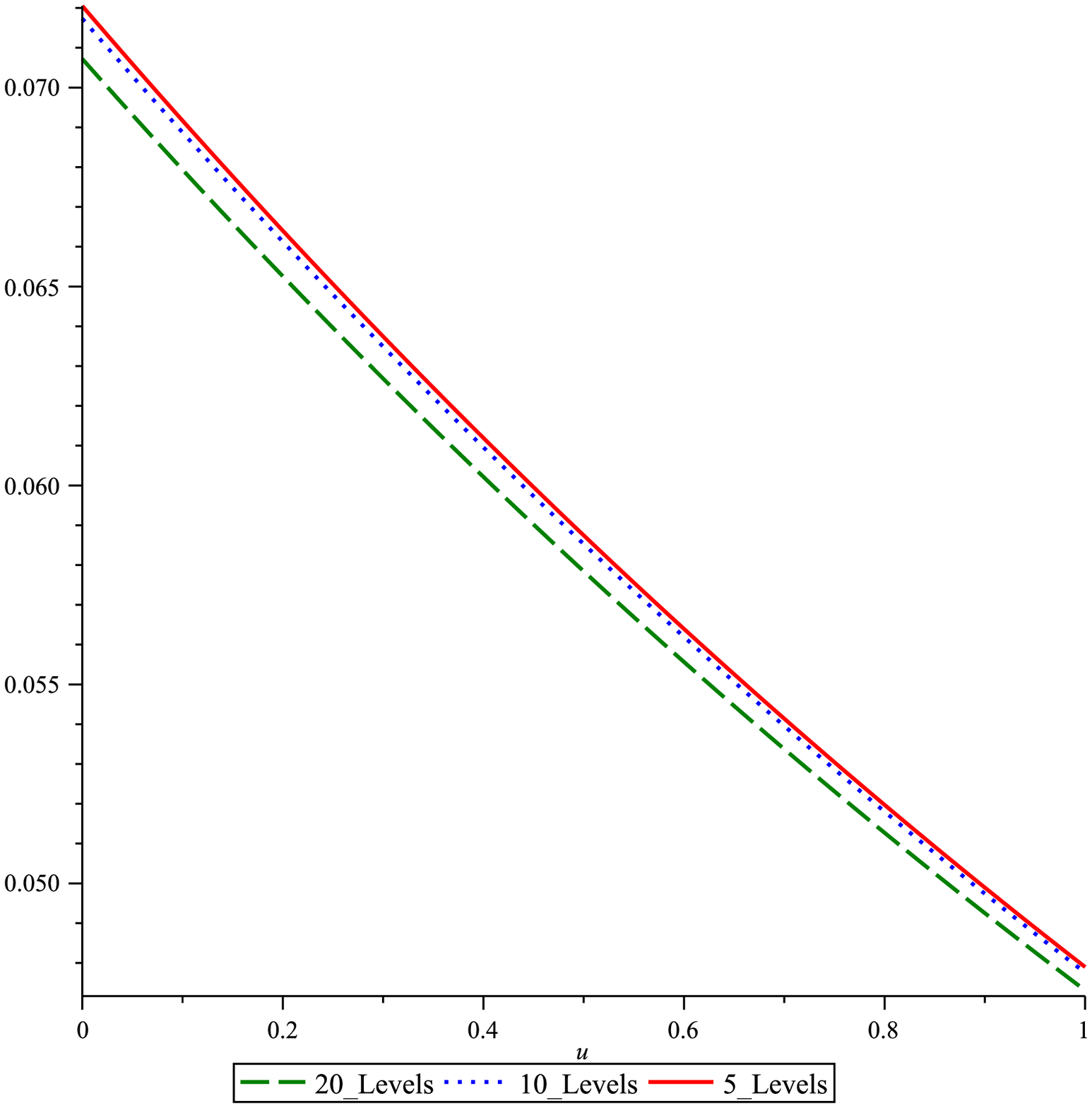}}\subfigure[]{
\includegraphics[width=4cm,height=5cm]{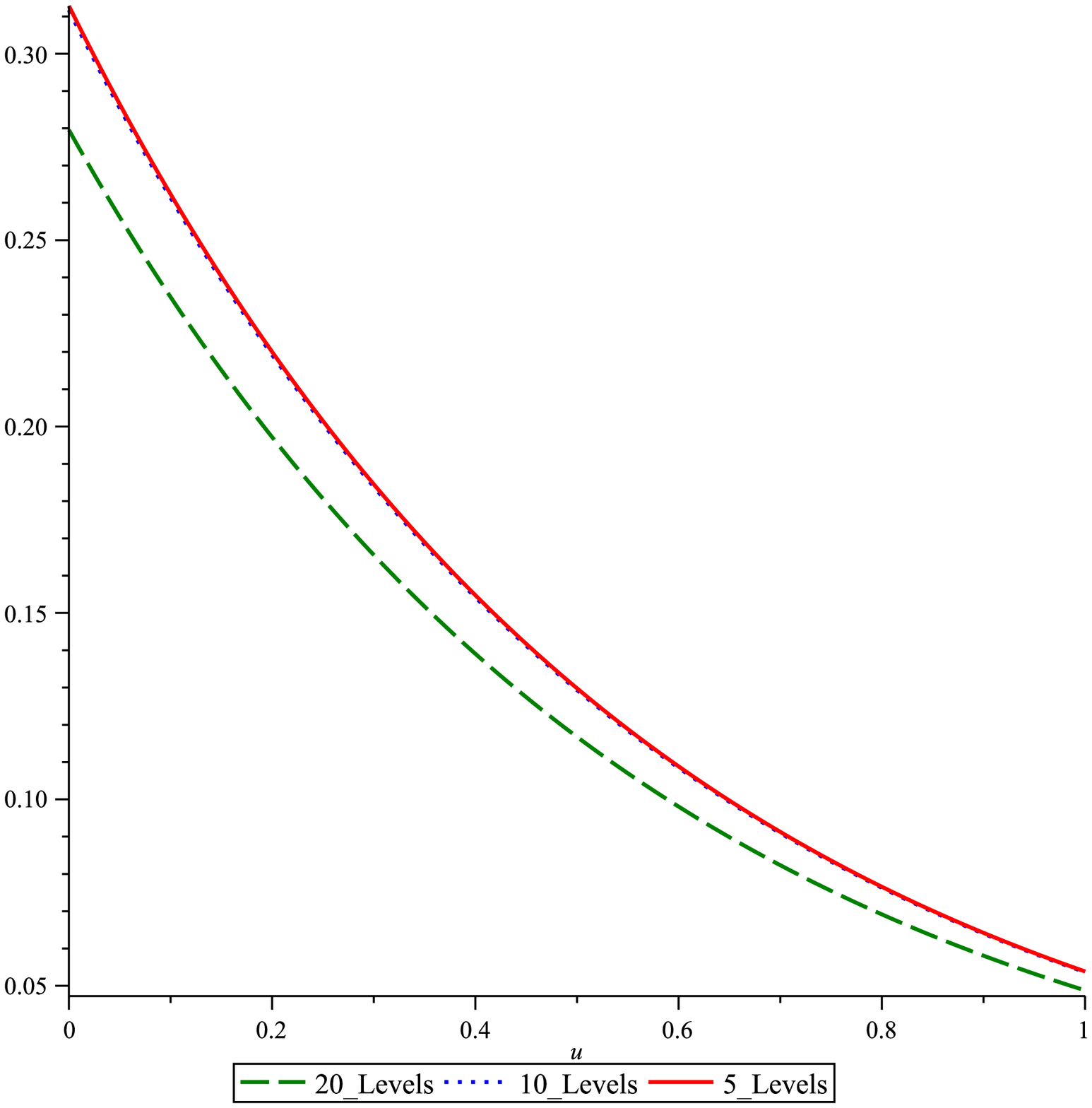}}\subfigure[]{
\includegraphics[width=4cm,height=5cm]{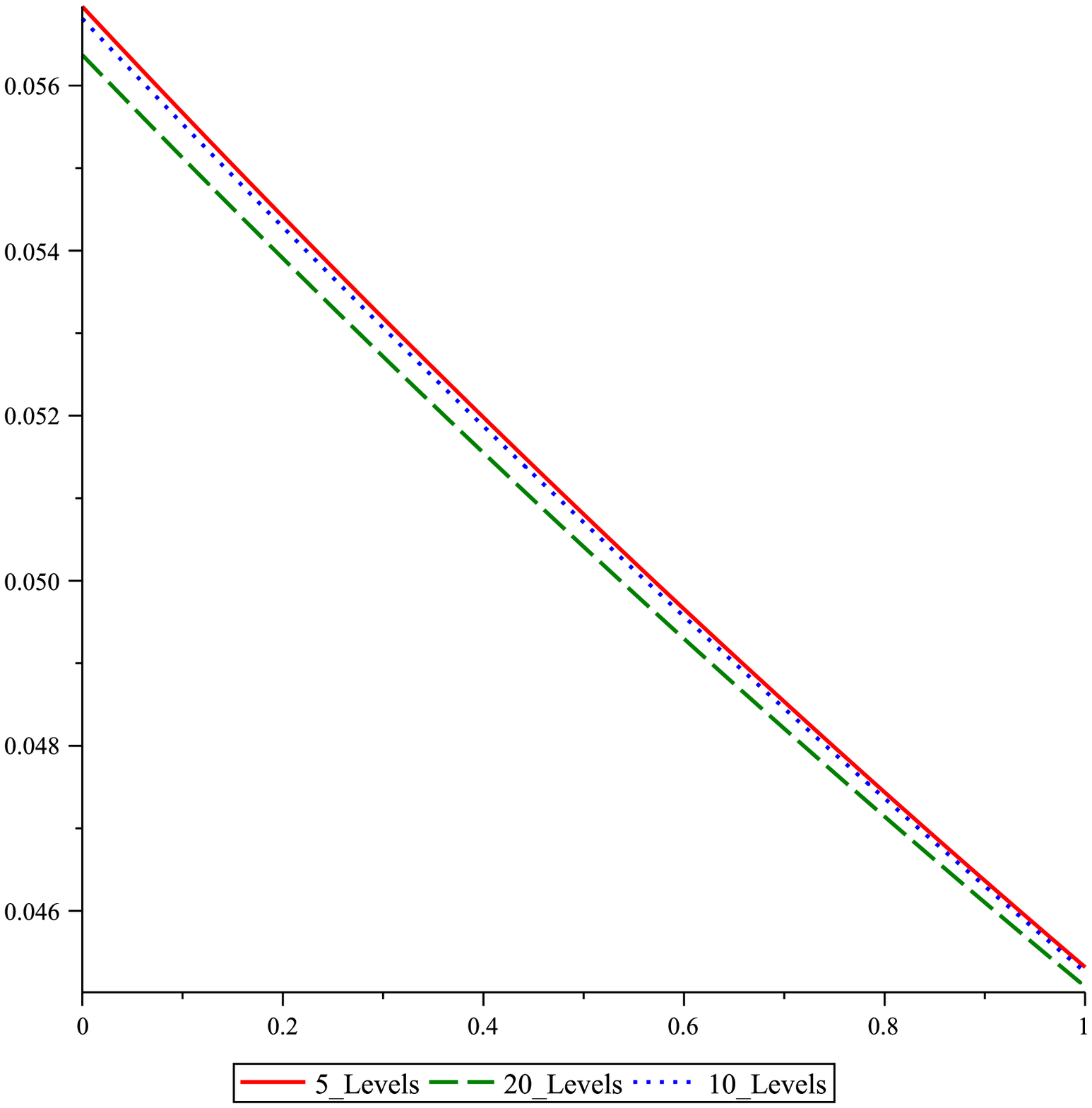}}\subfigure[]{
\includegraphics[width=4cm,height=5cm]{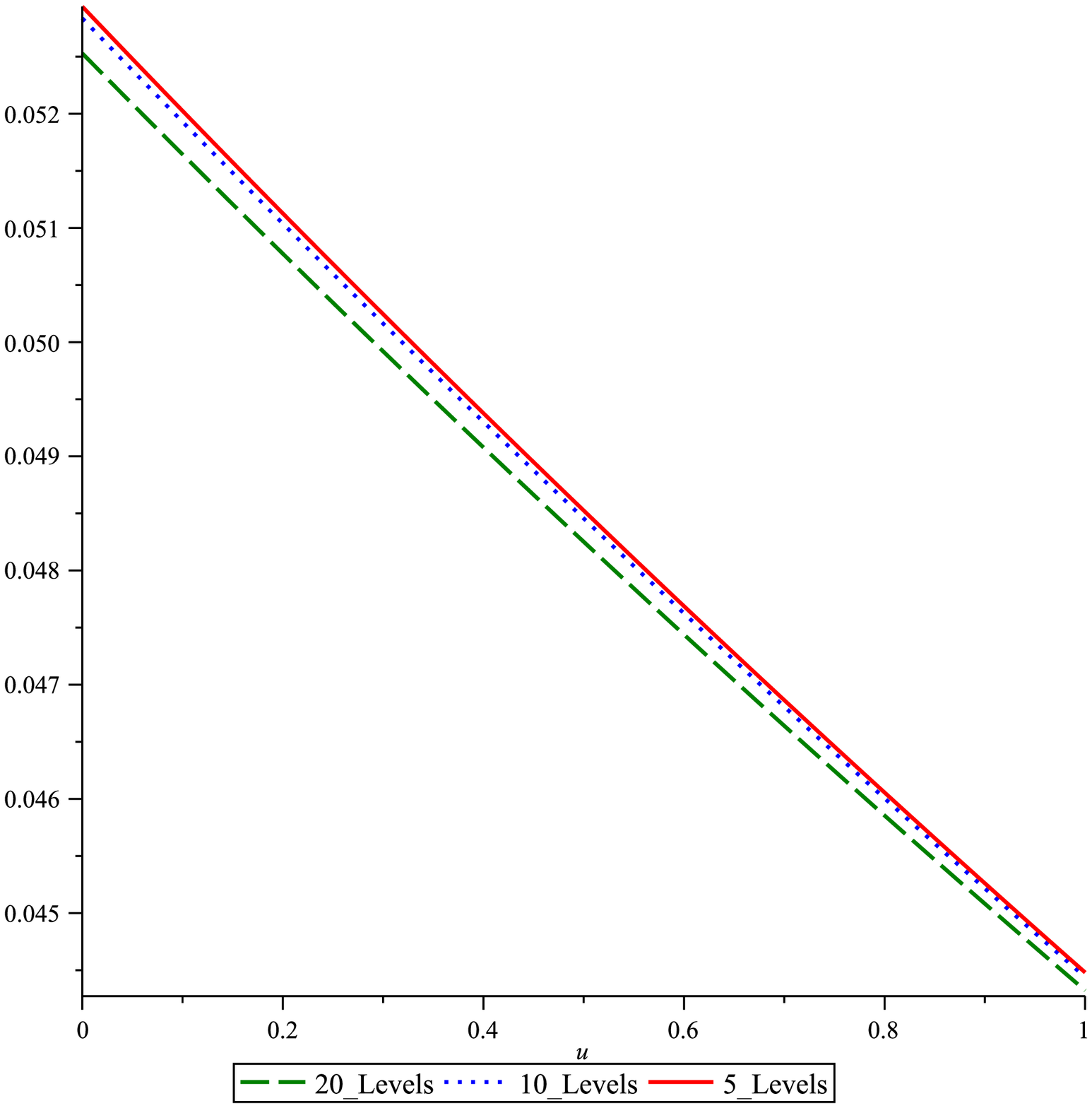}}
\caption{\scriptsize Figures (a) Density function of claim size
distributions; (b) to (d) comparison between the ruin probability
of such claim size distributions with respect to 5, 10, and 20
premium values scenarios, respectively; (e) to (h) comparison
between the ruin probability of 5, 10, and 20 premium values
scenarios with respect to different claim size distributions (i.e,
Gamma(1,1), Gamma(1,3), Gamma(3,2), Gamma(5,3)).}
\end{figure}
\end{center}
From the above figures, one may conclude that the ruin probability
$\psi(u)$ as function of initial wealth, $u,$ goes to zero and its
decay rate depends on tail of its corresponding claim size
distribution. Moreover, the ruin probability decreases as the
number of premium levels increases. The second observation in
different context has been pointed by several authors, see Denuit,
et al. (2007) for more detail.
\end{example}
\begin{example}
Suppose the claim size is a considerable heavy-tailed distribution
with density function
\begin{eqnarray}
\label{mixture_density}
 f_X(x) &=&
  0.150e^{-(x-4)^2/2}+\frac{0.075x}{(1+x^2)^2}+0.074e^{-3x}+\frac{0.224e^{-1/(2x)}}{x^{1.5}},
  x\geq0
\end{eqnarray}
and random premium $C$ has the Gamma distribution with parameters
$\kappa=\theta=3.$ Moreover, suppose that the number of sold
contracts, $N_1(t),$ and number of arrived claims, $N_2(t),$ are
two independent Poisson processes with intensity $\lambda_1=18$
and $\lambda_2=11.$ Using a mathematical software such as Maple,
one may show that an equation $D(s)=0$ has infinity many simple
non-positive roots, see Figure 2(b) for an illustration.

For simplicity, we just considered roots which their both real and
imaginary parts are fall within the interval $[-1,2].$ In this
subset of complex plane ${\Bbb C},$three roots
$z_1:=-0.5639909305-1.521665917I,$ $z_2:=-0.4233614172,$ and
$z_3:=-0.5639909305+1.521665918I$ have been found. Now using
Theorem \ref{linear_systems}, one may verify that the ruin
probability can be evaluated using
$\psi(u)=A_1e^{z_1u}+A_2e^{z_2u}+A_3e^{z_3u},$ where\footnotesize{
\begin{eqnarray*}
A_1&=&-\frac{-M_1z_3\beta_{23}z_2\beta_{32}+M_1T_2T_3+z_1\beta_{21}z_2\beta_{32}M_3+z_1\beta_{21}M_2
T_3+z_1\beta_{31}T_2M_3+z_1\beta_{31}z_3\beta_{23}M_2}{T_1z_3\beta_{23}z_2\beta_{32}-T_1T_2T_3+T_2z_3\beta_{13}
z_1\beta_{31}+z_2\beta_{21}^2z_1T_3+z_3\beta_{23}z_2\beta_{21}z_1\beta_{31}+z_3\beta_{13}z_1\beta_{21}z_2\beta_{32}}\\
A_2&=&\frac{-\beta_{21}\beta_{31}M_3z_2z_1-\beta_{21}z_2M_1T_3+\beta_{31}z_3\beta_{13}M_2z_1-z_2T_1\beta_{32}M_3-z_2M_1z_3
\beta_{13}\beta_{32}-T_1M_2T_3}{T_1z_3\beta_{23}z_2\beta_{32}-T_1T_2T_3+T_2z_3\beta_{13}z_1\beta_{31}+z_2\beta_{21}^
2z_1T_3+z_3\beta_{23}z_2\beta_{21}z_1\beta_{31}+z_3\beta_{13}z_1\beta_{21}z_2\beta_{32}}\\
A_3&=&\frac{-z_3\beta_{23}z_2\beta_{21}M_1+T_1T_2M_3+T_1z_3\beta_{23}M_2-z_2\beta_{21}^2z_1M_3+T_2z_3\beta_{13}M_1+z_
3\beta_{13}z_1\beta_{21}M_2)}{(T_1z_3\beta_{23}z_2\beta_{32}-T_1T_2T_3+T_2z_3\beta_{13}z_1\beta_{31}+z_2\beta_{21}^2
z_1T_3+z_3\beta_{23}z_2\beta_{21}z_1\beta_{31}+z_3\beta_{13}z_1\beta_{21}z_2\beta_{32}},
\end{eqnarray*}}\normalsize
where $T_j:=\lambda_2z_jM^\prime_X(-z_j)),$
$M_j:=\lambda_2M_X(-z_j)-\lambda_2,$ and
$\beta_{ij}=(\lambda_1+\lambda_2-\lambda_1M_C(z_i)-\lambda_2M_X(-z_j))/(z_j-z_i)$
for $i,j=1,2,3.$ Therefore, the ruin probability can be
approximated by
\begin{eqnarray*}
  \psi(u_0)
  &=&\frac{1}{100}||(-0.26+0.52I)e^{(-0.5640-1.5217I)u_0}-2.87e^{-0.4234u_0}+(-0.26-0.52I)e^{(-0.5640+1.5217I)u_0}||.
\end{eqnarray*}
The Cram\'er-Lundberg upper bound for this situation is
$e^{-0.4234u_0}$ which is significantly improved by our method, see
Figure 2(c).
\begin{center}
\begin{figure}[h!]
\centering\subfigure[]{
\includegraphics[width=5cm,height=5cm]{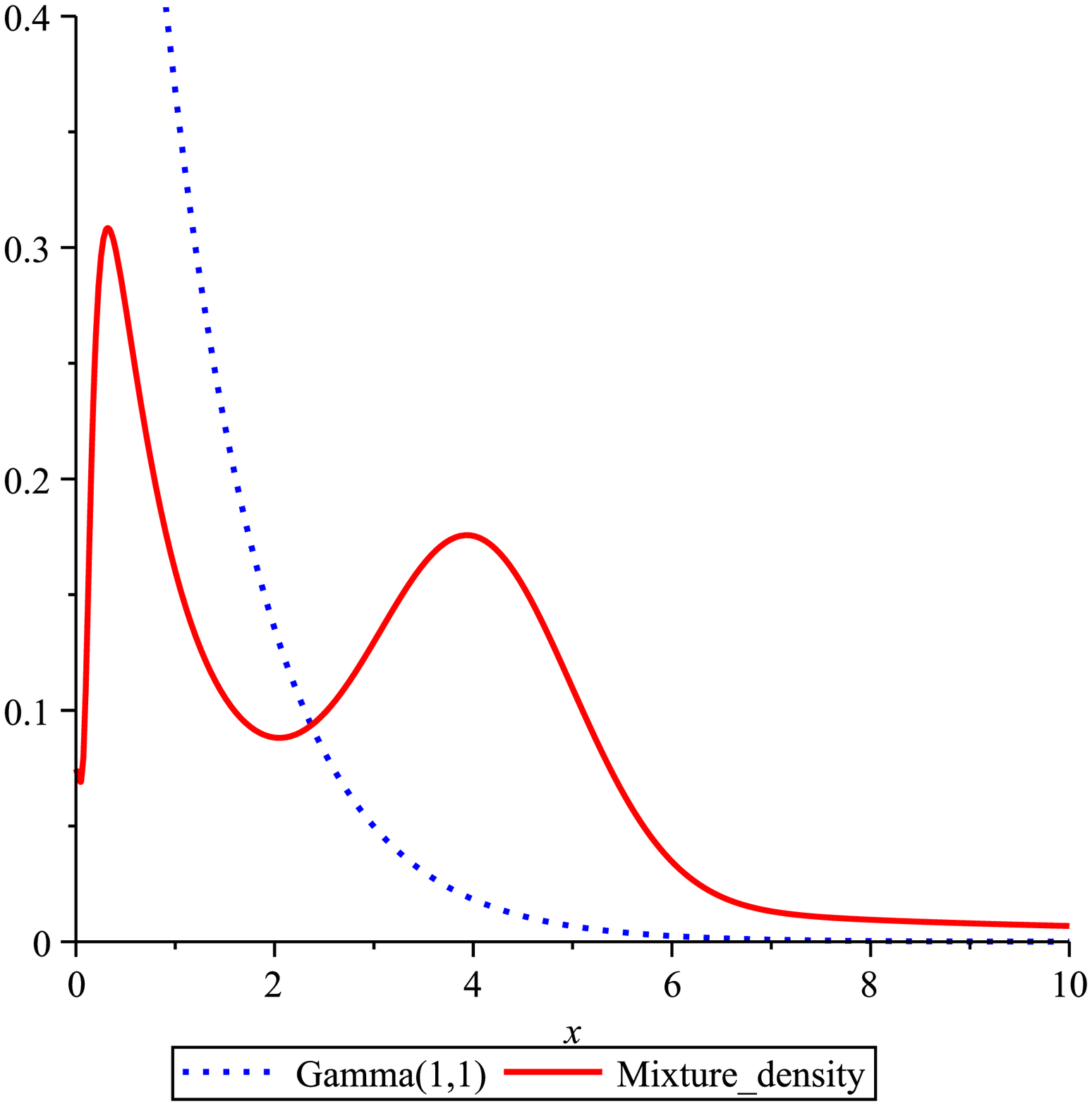}}\subfigure[]{
\includegraphics[width=5cm,height=5cm]{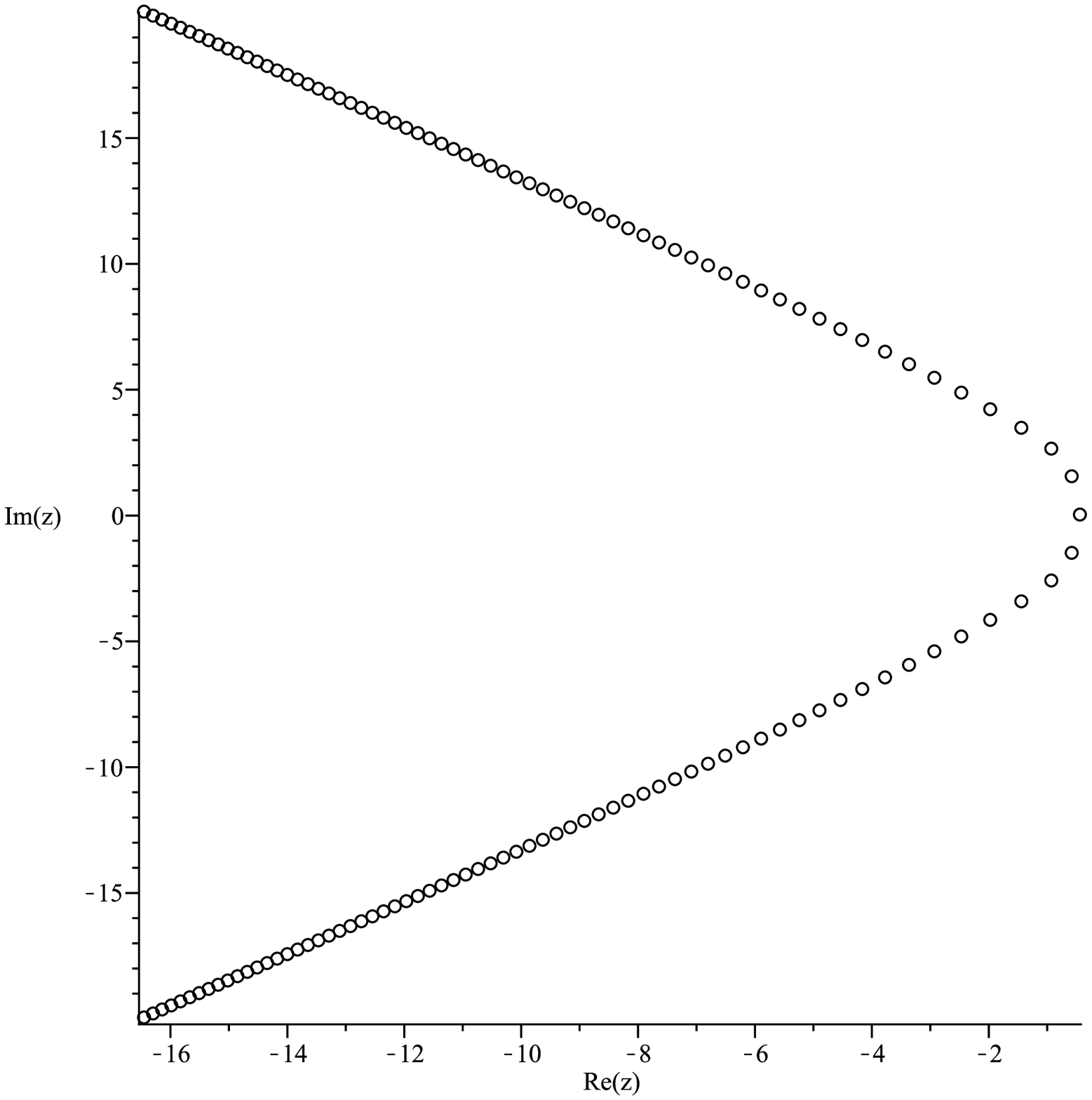}}\subfigure[]{
\includegraphics[width=5cm,height=5cm]{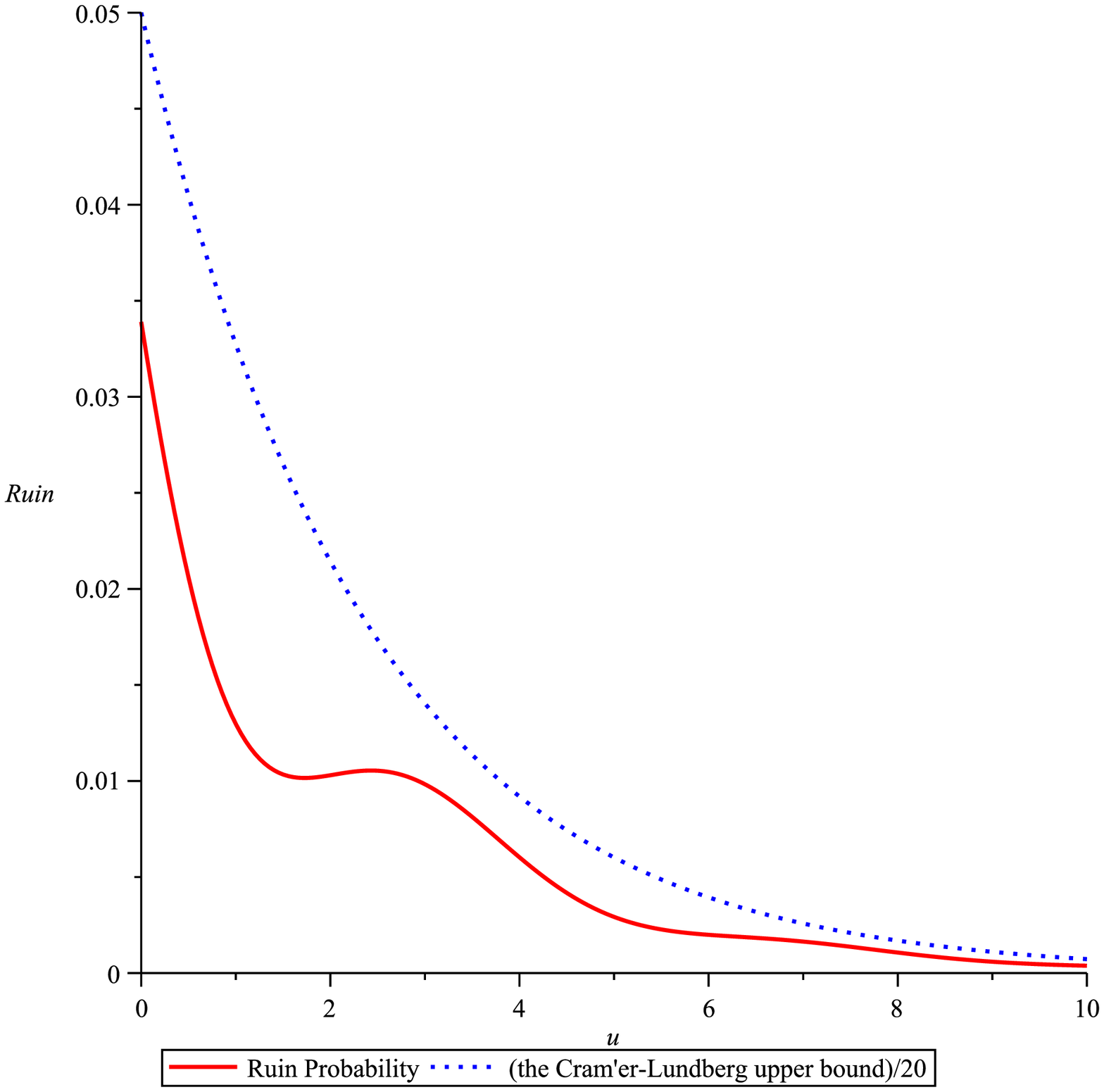}}
\caption{\scriptsize Figures (a) Density function of mixture claim
size density \ref{mixture_density} accompanied with Gamma (1,1)
density function; (b) Roots of equation $D(z)=0$ where $|z|\leq
20\sqrt{2};$ and (c) $\frac{1}{20}\times$ the ruin probability and
the Cram\'er-Lundberg upper bound of the surplus process
\ref{surplus-process-BMS} under mixture claim size density
\ref{mixture_density}.}
\end{figure}
\end{center}
\end{example}
\section{Conclusion and suggestions}
The Ruin theory provides some potential tools to study today's
solvency assessments for non-life insurance companies, see
W\"uthrich (2014) for a practical application. In most practical
applications the ruin probability of surplus process
\ref{surplus-process-BMS} has been approximated by the
Cram\'er-Lundberg upper bound $e^{-Ru_0},$ where adjustment
coefficient $R.$

This article assumed the ruin probability $\psi(\cdot)$ is an
exponential type $T$ and $L_2({\Bbb R})$ function. Then, it
derived several approximate formulas for the ruin probability of a
double stochastic compound Poisson process. The approximated ruin
probability constructed based upon roots of equation $D(s)=0.$
Certainly, $-R$ is one of such roots. Therefore, the
Cram\'er-Lundberg upper bound, multiplied by a constant, will be a
part of our approximated ruin probability. Since $-R$ is the
largest roots of the equation $D(s)=0,$ our approximation method
improved the Cram\'er-Lundberg upper bound. Such improvement is
significant whenever equation $D(s)=0$ has more than one simple
non-positive root, see Example 2.

The exponential type assumption can be dropped and a formal
solution can be verified by substitution into Equation
\ref{integro-differential-equation}. But, the generally realistic
$L_2({\Bbb R})$ assumption cannot be dropped, because our methods
always produce solutions with finite $L_2({\Bbb R})$ norm.
Furthermore, it is essential to our methods that the Laplace
transform of $f_X$ does not have any branch points in the complex
plane.
\section*{Acknowledgements}
The support of Shahid Beheshti University and Natural Sciences and
Engineering Research Council (NSERC) of Canada are gratefully
acknowledged by authors.

%Many thanks goes to an anonymous reviewer who improved this
%article by his/her constructive comments.

\end{document}